\documentclass[twocolumn]{revtex4}
\usepackage{latexsym}
\usepackage{amsmath,amsthm,amssymb,amsfonts,amsbsy,amscd,amstext}
\usepackage{times}
\usepackage{amssymb}
\usepackage{fancyheadings}
\usepackage[T1]{fontenc}
\usepackage[utf8]{inputenc}
\usepackage{fancyhdr}
\usepackage{url}
\usepackage{hyperref}
\usepackage{color}
\usepackage{verbatim}
\usepackage{algorithm,algorithmic}
\usepackage{pgf,tikz}
\usepackage{enumerate}
\usepackage{cases}
\usepackage{verbatim}
\usepackage{bm}
\usepackage{textcomp}
\usepackage{graphicx}
\usepackage{xcolor}

 \newtheorem{theorem}{Theorem}[section]

 \newtheorem{remark}[theorem]{Remark}
 
\newcommand{\J}{\mathbf{J}_e}
\newcommand{\DZ}{\Delta\mathbf{Z}}
\newcommand{\Z}{\mathbf{Z}}
\newcommand{\Om}{\Omega}
\newcommand{\bnu}{\boldsymbol{\vec{n}}}

\newcommand{\p}{\cdot}
\newcommand{\dv}{\,\delta v}
\newcommand{\ds}{\,\delta s}

\newcommand{\A}{\mathbf{A}}

\newcommand{\V}{\nabla\textbf{V}_c}
\newcommand{\Vv}{\textbf{V}_c}
\newcommand{\Vh}{\textbf{V}_{c,\natural}}
\newcommand{\As}{\mathbf{A}_{\tau}}
\newcommand{\PhiS}{\Phi_{\tau}}
\newcommand{\Vs}{\nabla_{\tau}\textbf{V}_c}

\newcommand{\E}{\mathbf{E}}
\renewcommand{\H}{\mathbf{H}}
\newcommand{\curl}{\textbf{curl\hspace{.01in}}}
\newcommand{\dvg}{\textbf{div}}

\newcommand{\grad}{\nabla}
\newcommand{\gradS}{\nabla_{\tau}}

\usepackage{ff++listings}
\begin{document}


\title{\bf A Fast Eddy-current Non Destructive Testing Finite Element Solver in Steam Generator}
\author{Mohamed Kamel RIAHI$^{\star}$}\thanks{Corresponding author. Email: riahi@njit.edu,  Tel.: +1 (973)-5966-084.}
\affiliation{$^\star$Department of Mathematical Science, New Jersey Institute of Technology, Newark, New Jersey, USA}
\date{\today}

 
\begin{abstract}
\noindent 
In this paper we present an advanced numerical method to simulate a real life challenging industrial problem that consists of the non-destructive testing in steam generators. We develop a finite element technique that handles the big data numerical set of systems arising when a discretization of the eddy-current equation in three dimensional space is made. Using a high performance technique, our method becomes fully efficient. We provide numerical simulations using the software Freefem++ which has a powerful tool to handle finite element method and parallel computing. We show that our technique speeds up the simulation with a good efficiency factor. 

\vspace*{2ex}\noindent\textit{\bf Keywords}: Maxwell's Equation, Non-destructive testing, Eddy-Current approximation, Numerical Methods, Parallel Algorithm, High performance computing, Industrial problem.
\\[3pt]
\\[3pt]
\end{abstract}

\maketitle

\thispagestyle{fancy}

\section{Introduction and industrial model}\label{Mission of the Journal}
\vspace*{0.2cm}

Many industrial applications use Eddy-Current Testing (ECT) to evaluate cracks, defaults and other types of anomalies in the material engines. In general, ECT is a non-invasive technique where electric impedance measurement plays an important role in the evaluation of the quality of the tested material. Eddy currents are created through a process called electromagnetic induction. When alternating current is applied near the conductor (such as copper wire), a magnetic field develops in and around this conductor. The measurement of the electrical impedance enables a detection of anomalies in the tested pieces.    
The computer simulation of the eddy-current problem has an enormous role in the automatization of the ECT, especially when inverse type problems are conducted, such as shape identification of deposits and cracks detection~\cite{1386227,RPQNE,bendjoudi,HuangTakagiFukutomi,HuangTakagi,girarclogging,haddar2015axisymmetric} (we may also refer to~\cite{ndtdatabase} for further lecture). To do so, the direct problem needs to be robust and effective. 

	In this work, we are concerned with the numerical simulation on parallel computers of the direct problem arising from the ECT of a steam generator (SG). The numerical problem involves the finite element approximation of the eddy-current equations in a time harmonic low frequency regime. Because of the vectorial aspect of the three dimensional unknown, the discretization of the problem leads to a huge and ill-conditioned system. Its numerical resolution is time and memory consuming. This has motivated us to construct an algorithm based on high performance computing tools.
	
	SGs are critical components in nuclear power plants. Heat produced in a nuclear reactor core is transferred as pressurized water at high temperature via the primary coolant loop into an SG, consisting of tubes in U-shape, and boils coolant water in the secondary circuit on the shell side of the tubes into steam. The steam is then delivered to the turbine generating electrical power. The SG tubes are held by the broached quatrefoil tube support plates (TSP) with flow paths between tubes and plates for the coolant circuit.

	Without disassembling the SG, the lower part of the tubes -- which is very long -- is inaccessible for normal inspections. Therefore, an ECT procedure is widely practiced in industry to detect the presence of defects, such as cracks, flaws, inclusions and deposits.

	In this paper, we go over a specific technique and tools using high performance computing. We provide an efficient and scalable finite element algorithm to solve the ECT problem in the SG. We consider a potential formulation supplemented by a Coulomb gauge condition, where the electric field is represented in a suitable space by a magnetic vector potential and scalar electric potential (see for instance~\cite{MR2680968} and reference therein). This involves a coupled system between the new variables, which we solve using a sophisticated technique. 
	
	In the SG case, the ECT procedure consists in withdrawing a probe constituted by two coils that move together all the way through the tube. These coils are called generator and receiver coils, where the generator coil produces the source term current excitation that produces a magnetic field, which in turn penetrate materials nearby and produce an eddy-current by induction. The receiver coil, has a detective role, where it can measure the variation of electric impedance due to the induced eddy-current. The motivation behind our parallelism technique is that the ECT procedure is a complicated task that needs a careful attention starting from the mesh generation regarding the motions of the coils along the scan zone, to the resolution of the coupled system that governs the electric field solution of the eddy-current problem. Indeed, we consider the vectorial $\curl\curl$ version of the eddy-current equations in three dimension. Thus the finite element discretization of the problem at hand may significantly lead to a huge linear coupled system, especially when mesh refinement is needed in order to have a good approximation of the exponential decay of the wave when it penetrates materials. 

Our work is presented as follows.  
In section 2 we present the eddy-current equations. We reformulate the problem using potential formulations, then we give the system that has to be solved. Our presentation uses bi linear formulation obtained after setup of variational formulation of the eddy-current problem. In section 3 we  
In section 3 we describe the numerical method that we use to deal with the particularity of the ECT in SG, where a set of operations has to be done and that involves numerical solution of eddy-current problem. We provide sample of Freefem++  scripting that show how to implements mathematical formulas using the aforementioned software. We give numerical results in section4, that shows the efficiency of our approach in the speedup of the ECT. 
\vspace*{0.2cm}
\noindent

\section{Time harmonic Eddy-current equations}
\vspace*{0.2cm}
  \noindent
  
  The following description of the eddy-current problem follows \cite{MR2680968}.
  
  The time harmonic Eddy-current approximation of Maxwell's equations for the electric field $\E$ and the magnetic field $\H$ and a source excitation term $\J$ reads:
\begin{equation}\label{MAMF}
\begin{cases}
\curl\H-\sigma \E=\J&\text{ on } \Omega,\\
\curl\E - i\omega\mu \H = 0.&\text{ on } \Omega.
\end{cases}
\end{equation}
where $\sigma$ stands for the electric conductivity, $\mu$ stands for the magnetic permeability and $\omega$ stands for the frequency. The computational domain $\Om\subset\mathbb{R}^3$ is a bounded with Lipchitz boundary $\partial\Om$. The above system represents Maxwell-Amp\`ere and Maxwell-Faraday equations, where the current displacement term has been dropped to arrive at the eddy-current model~\cite{Ammari:2000:JEC:354423.354457,bossavit2004electromagnetisme}. In addition, in order to insure uniqueness of the solution; Eq.\eqref{MAMF} is thus supplemented with a boundary condition. In our case we consider a perfect magnetic boundary condition i.e. $\H\times\bnu=0$ on $\partial\Om$ with $\bnu$ stands for the unit normal vector pointing outward at the domain $\Om$. Our computational domain is decomposed, with respect to the conductivity $\sigma$, to a conductor material $\Om_C$ characterized with $\sigma\neq0$ and insulator $\Om_I$. Thus $\Om=\Om_I\cup\Om_C$, where the common interface is denoted as $\Gamma=\Om_I\cap\overline{\Om_C}$. Remark that $\Gamma\cap\partial\Om\neq\{0\}$. 

From Eq.~\eqref{MAMF}$_1$ one can see that $\curl\H = \J$ in the insulating region $\Om_I$. This imposes a condition on the excitation term $\J$ where it is required to be a divergence free vector field. 
 
  
	In order to cope with the geometry configuration of the computational domain, it is common to use potential formulation to better handle the eddy current problem numerically. 
In the sequel we consider the magnetic vector potential formulation by introducing the magnetic vector potential $\A$ and the scalar electric potential $\Vv$ (uniquely defined on the conductor $\Om_c$)~\cite{biro2007coulomb} where they satisfy 
\begin{equation}\label{defA}
\begin{cases}
\E = i\omega\A +  \V &\text{ on } \Om,\\
\mu\H= \curl \A &\text{ on } \Om.
\end{cases}
\end{equation}
In order to avoid singularity in this system, regarding the change of variables Eq.\eqref{defA}, and ensure well-posedness; in the sense of the magnetic potential vector $\A$ is unique, it is classical and necessary to impose additional gauge conditions. In this work we consider Coulomb gauge conditions which read
\begin{equation}\label{CGCdiv}\dvg \A=0\text{ in } \Om,
\end{equation}
We also need to impose the boundary condition $\A\p\bnu = 0$ on $\partial\Om$ in order to close the problem. A classical technique~\cite{MR2298698} incorporates the constraint Eq.\eqref{CGCdiv} using a penalization term $-\frac 1{\tilde\mu}\grad\dvg\A$, in the Amp\`ere equation, where $\tilde\mu$ is a suitable average of $\mu$ in $\Om$. Therefore a strong formulation ($\A,\V$) of our eddy-current problem writes
 \begin{equation}\label{strongPbAV}
 \begin{cases}
 \curl\big( \dfrac 1 \mu\curl \A\big) -\frac 1 {\tilde\mu}\grad\dvg\A-\sigma (i\omega\A - \V) = \J \!\!\!\!&\text{in }\Om,\\
\dvg \big( i\omega\sigma\A+\sigma\V \big) = 0&\text{in }\Om_c,\\
\big( \sigma i\omega\A+\sigma\V\big)\p\bnu = 0&\text{on }\Gamma,\\
\A\p\bnu = 0 &\text{on }\partial\Om,\\
\big(\dfrac{1}{\mu}\curl\A\big)\times\bnu = 0 &\text{on }\partial\Om,
 \end{cases}
 \end{equation}
where $\Vv$ is determined up to an additive constant.

	Consider the space $H(\textbf{curl},\Omega)\cap H_0(\textbf{div},\Omega)$.
where $H(\textbf{curl};\Omega):=\{{\bf u}\in (L^2(\Omega)\big)^3 \,|\, \curl {\bf u}\in(L^2(\Omega)\big)^3 \},$ and $H(\textbf{div};\Omega):=\{ {\bf u}\in(L^2(\Omega)\big)^3\,|\, \nabla\p {\bf u}\in L^2(\Omega) \},$
also we have $ H_0(\textbf{div};\Omega):=\{{\bf u}\in H(\textbf{div};\Omega)\, |\, {\bf u}\p\bnu_{|\partial\Omega}=0\}.$

	Let us take test functions $\Phi\in H(\textbf{curl};\Omega)\cap H_0(\textbf{div};\Omega)$ and $\varphi\in H^1(\Omega_c)$ for the Eq.~\eqref{strongPbAV}$_1$ and the Eq.~\eqref{strongPbAV}$_2$ respectively. After integrating by part we obtain the following weak formulations:
\begin{eqnarray}
\int_{\Omega}\dfrac{1}{\mu}\curl\A\p\curl\overline{\Phi} \dv +\dfrac{1}{\tilde\mu}\int_{\Omega}\dvg\A\dvg\overline{\Phi} \dv \notag\\
	\hspace{-2cm}- \int_{\Omega_c}\sigma(i\omega\A+\V)\p\overline{\Phi} \dv = \int_{\Omega}\J\p\overline{\Phi} \dv\label{varfcurl} \\
	\hspace{-2cm}\int_{\Omega_c}\sigma\big(i\omega\A+\V\big)\p\overline{\nabla\varphi} \dv = 0\label{varflaplace} .
\end{eqnarray}

	We multiply Eq.~\eqref{varflaplace} by $\dfrac{-1}{i\omega}$ to obtain : 
\begin{equation*}
\displaystyle\dfrac{-1}{i\omega}\int_{\Omega_c}\sigma\big(i\omega\A+\V\big)\p\overline{\nabla\varphi} \dv= 0.
\end{equation*}	
and couple this with Eq.~\eqref{varfcurl} in a single mixed weak variational formulation, which is written:
\begin{align*}
&\int_{\Omega}\dfrac{1}{\mu}\curl\A\p\overline{\curl\Phi} \dv +\dfrac{1}{\tilde\mu}\int_{\Omega}\dvg\A\overline{\dvg\Phi} \dv\\
	&- \dfrac{1}{i\omega}\int_{\Omega_c}\sigma(i\omega\A+\V)\p (i\omega\overline\Phi+\overline{\grad\varphi}) \dv \\
&=\int_{\Omega}\!\!\!\J\p\overline{\Phi} \dv.
\end{align*}
	For sake of simplicity, we define the sesquilinear form $\mathcal{L}(\A,\Vv,\Phi,\varphi)$ as the right-hand side of the above, which is written:
\begin{eqnarray}\label{sesquil}\displaystyle
\mathcal{L}\big(\A,\Vv,\Phi,\varphi\big):=\int_{\Omega}\dfrac{1}{\mu}\curl\A\p\overline{\curl\Phi} \dv \\+\dfrac{1}{\tilde\mu}\int_{\Omega}\dvg\A\overline{\dvg\Phi} \dv \notag\\
	- \dfrac{1}{i\omega}\int_{\Omega_c}\sigma(i\omega\A+\V)\p (i\omega\overline{\Phi}+\overline{\grad\varphi}) \dv \notag.
\end{eqnarray}

and which represents a coupled system between the two potential variables. We quote hereafter the part that constitutes the coupled system given as variational forms

\begin{eqnarray} 
\mathcal{L}_{11}\big(\A,\Phi\big)=\int_{\Omega}\dfrac{1}{\mu}\curl\A\p\overline{\curl\Phi} \dv,\\+\dfrac{1}{\tilde\mu}\int_{\Omega}\dvg\A\overline{\dvg\Phi} \dv \notag\\
\mathcal{L}_{12}\big(\Vv,\Phi\big)=-\int_{\Omega_c}\sigma \V\p\overline{\Phi} \dv.\\
\mathcal{L}_{21}\big(\A,\varphi\big)=-\int_{\Omega_c}\sigma\A\p \overline{\grad\varphi} \dv,\\
\mathcal{L}_{22}\big(\Vv,\varphi\big)=- \dfrac{1}{i\omega}\int_{\Omega_c}\sigma\V \p\overline{\grad\varphi} \dv.
\end{eqnarray}

	 We describe briefly in the sequel how one can dismiss the volume part of the TSP by using the appropriate boundary condition since it has high conductivity as compared with the tube, its corresponding skin depth is then very small.
	 
\begin{figure}[!htpb]
\includegraphics[height=8cm,width=4cm]{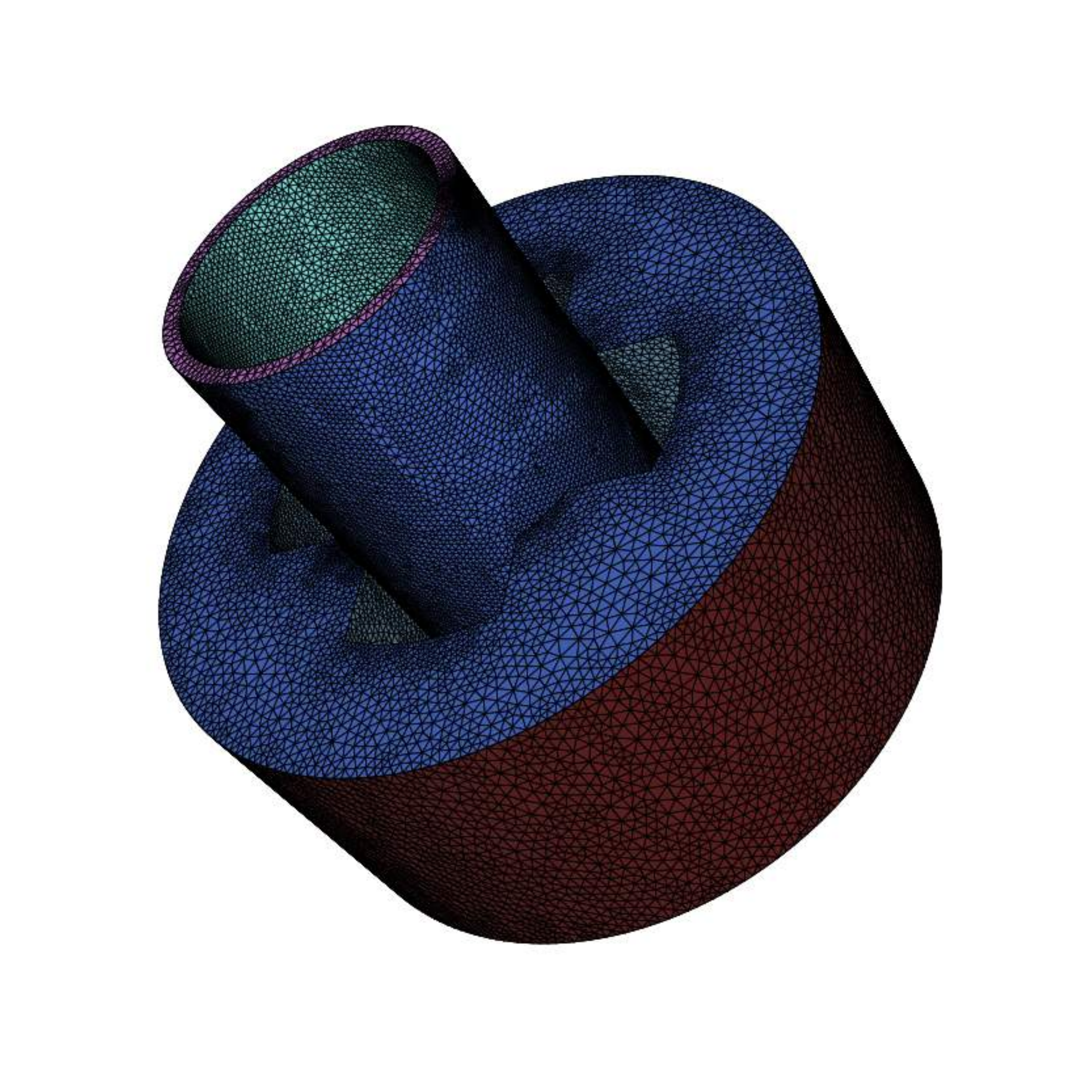}
\includegraphics[height=8cm,width=4cm]{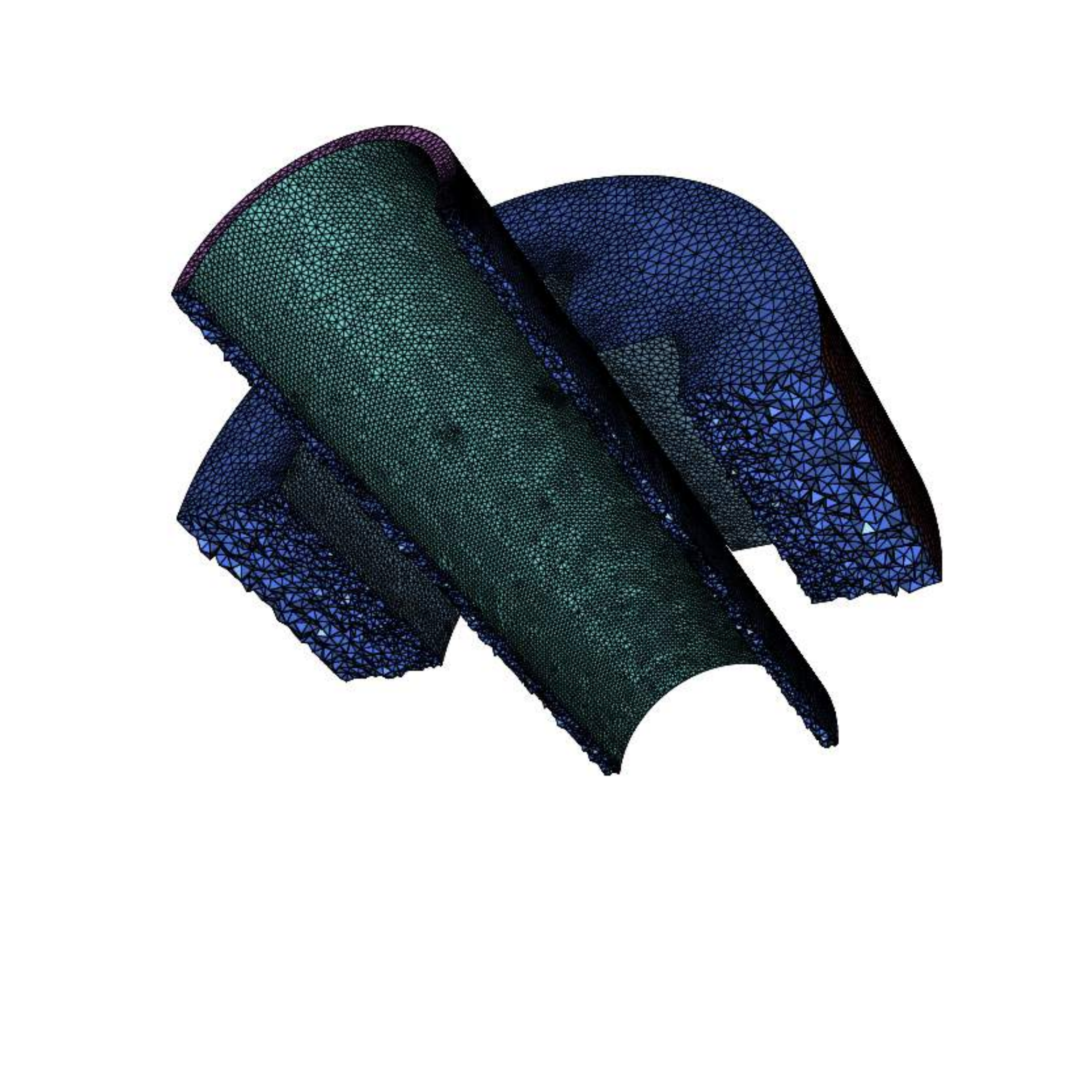}
\caption{Triangulation of the conductive materials in the steam generator. Tube and Tube-support-plate (left), and their respective vertical-slice (right).}\label{condpart}
\end{figure}

	Taking into account the effect of the TSP using the 3D model, described above, which requires a very thin mesh size (proportional to the skin depth) inside the TSP and leads to a huge size of the discrete 3D problem. We hereafter explain how one can avoid integration over the volume of the TSP by imposing an appropriate impedance boundary condition (IBC) on its boundary $\Gamma_p$. More precisely it is shown in~\cite{durufle2006higher} that electromagnetic field satisfies 
 \begin{equation}\label{IBC}
\bnu\times\frac{1}{\mu}\curl\A = -\dfrac{1}{\mathcal{Z}_{\Gamma_{p}}}(i\omega\A_\tau+{\bf \grad V}_{c,\tau}),  \quad\text{ on } \Gamma_{p}.
\end{equation}
(up to $O(\delta^2)$) on $\Gamma_{p}$. In the above equation; $\mathcal{Z}_{\Gamma_{p}}:= (1-i)\slash(\delta\sigma)$ with the skin depth $\delta:=\sqrt{ 2\slash(\omega\mu\sigma)}$ and the tangential component of the vector field $\A_T=\bnu\times\big( \A\times\bnu\big)$ (same apply for $\Vv$).
 Therefore, if $\delta$ is sufficiently small i.e. $\omega\sigma_p\mu$ is sufficiently large and Eq.~\eqref{IBC} is a very good approximation.

	The impedance surface term of the weak formulation of Eq.~\eqref{strongPbAV} at the interface $\Gamma_p$ of the TSP is written:
\begin{equation}\label{ibcROTA}\footnotesize
\int_{\Gamma_{p}}\!\!\!\!\!\!\! \big(\bnu\times(\frac{1}{\mu_p}\curl\A)\big)\p\overline{\PhiS} \ds=\dfrac{-1}{\mathcal{Z}_{\Gamma_{p}}}\!\!\! \int_{\Gamma_{p}}\!\!\!\! (i\omega\As+\Vs)\p\overline{\PhiS} \ds
\end{equation}
\begin{equation}\label{ibcNORM}\footnotesize
\int_{\Gamma_{p}}\!\!\!\!\!\!\! \sigma_p(i\omega\A+\V)\p\bnu\,\overline{\varphi} \ds =\\\dfrac{-1}{\mathcal{Z}_{\Gamma_{p}}}\int_{\Gamma_{p}}\!\!\!\!(i\omega\As+\Vs)\p\overline{\gradS\varphi} \ds,
\end{equation}
We refers for instance to~\cite{haddar:hal-01044648} for details related to the update on the variational formulation with respect to the impedance boundary condition. Here, we shall consider the general case in the numerical experiments, which consist of the volume representation of the TSP in the variational formulation. 

The mathematical formulation for the evaluation of the electric impedance $\Z$ see~\cite{auld1999review}. It is common to use the following types of signals 
\begin{equation}\label{impedmod}
\begin{cases}
{\bf Z}_{FA} = \frac{i}{2}\big(\DZ_{11}+\DZ_{12}\big),\\
{\bf Z}_{F3} =\frac{i}{2} (\DZ_{11}-\DZ_{22}\big).
\end{cases}
\end{equation}
which are respectively the absolute signal mode and the differential  signal mode. The term $\Z_{kl}$ with $k,l$ in $\{1,2\}$ represents the volume impedance measured with the coil $k$ in the electromagnetic field induced by the coil $l$. It is written:
\begin{eqnarray}\label{impedkl}
\displaystyle
\DZ_{kl}:=\frac{1}{i\omega}\frac{\mu_\epsilon-\mu_d}{\mu_d\mu_\epsilon}\int_{\Om_{d}}\big( \curl\A_k\p\curl\A_l^{\epsilon}\big)\dv \notag\\
+ (\sigma-\sigma_{\epsilon})\!\!\!\int_{\Om_d}\!\!\!\!\!(i\omega\A_k + \V)\p(i\omega\A_{\tau,l}^{\epsilon}+{\bf \grad V}^{\epsilon}_{c,l,\tau}) \dv,
\end{eqnarray}
where $i\omega\A_{\tau,l}^{\epsilon}+{\bf \grad V}^{\epsilon}_{c,l,\tau}:=:\E_{l,\tau}^0$ refers to the tangential component of the electric field propagating in "vacuum". Rigorously one has to take care of the presence of loop field~\cite{MR3090156} in the vacuum. Their presence is related to the geometry of the computational domain (see~\cite{MR2680968,bossavit2004electromagnetisme} and reference therein).  
\begin{remark}
As eddy-current approximation of the Maxwell equation can be seen as the limit when the permittivity $\varepsilon$ tends to zero in Maxwell equations. We may use the vacuum as a medium with very low permittivity. In our case we consider a very low conductivity $\sigma_\epsilon\approx \epsilon<<1$ to describe the vacuum instead. 
\end{remark}

\section{\normalsize Numerical method}
We describe in this section our numerical method which consists of the resolution of the coupled systems that solve the eddy-current problems necessary for the evaluation of the electric impedance $\Z$ defined at Eq.~\eqref{impedmod}. We describe the major steps for the non-destructive test procedure and give details related to the resolution where numerical decisions have been made in order to better handle big ill-conditioned systems.

We suppose from now on, that the computational domain $\Omega$ is a Lipschitz polyhedra in $\mathbb{R}^3$, practically $\Omega$ represent the cylinder that envelops the Tube and the TSP. Once the computational domain is fixed, we then introduce a family of triangulation $\mathcal{T}_{\natural}$ of $\Omega$, the subscript $\natural$ stands for the largest length of the edges of the triangles that constitute $\mathcal T_{\natural}$. The Tetrahedrons of $\mathcal T_{\natural}$ match on the interface between the conductive part, i.e. tube and TSP ($\sigma\neq 0$) and the insulator part ($\sigma=0$). The triangulations of the conductive part are is in Figure~\ref{condpart}.

	 The electric field $\E\equiv(\A,\Vv)$ solution of the eddy current equation belongs to $H(\curl,\Omega) \cap H(\dvg,\Omega)\times H^{1}(\Omega_c)$. The fact that our computational domain (Cylinder as mentioned before) is a convex polyhedron is very important. It turns out that the space of smooth tangential vector fields $H_\tau^1(\Omega):=\big(H^{1}(\Omega)\big)^3\cap H^{}(\dvg,\Omega)$ coincides with the proper close subspace of $H(\curl,\Omega) \cap H(\dvg,\Omega)$. Thanks to this nice property, a finite element numerical approximation based on a nodal finite element will be considered for the electric vector potential $\Vv$ as well as for the magnetic vector potential $\A$~\cite{biro2007coulomb}. The infinite dimensional space $H_\tau^1(\Omega)\times H^{1}(\Omega)$ is therefore approximated with the finite-dimensional space $\mathcal W_\natural:=\mathcal V_\natural\times \mathcal C_\natural$, where $\mathcal V_\natural$ and $\mathcal C_\natural$ represent the discrete spaces of Lagrange nodal elements defined respectively as 
\begin{eqnarray*}
 \mathcal{V}_{\natural}:=\{ {\bf v}_\natural\in (C^0(\Omega))^3 | {\bf v}_{\natural_{|K}}\in (\mathbb P_1)^3, \forall K\in\mathcal T_\natural,\\ {\bf v}_\natural . \bnu=0 \text{ on } \partial \Omega\}
\end{eqnarray*}
 and $\mathcal{C}_{\natural}:=\{ c_\natural\in C^0(\Omega_c)| v_{\natural_{|K}}\in \mathbb P_1, \forall K\in\mathcal T_{c,\natural}\}$.
	
	 In the above $\mathbb P_1$ stands for the space of polynomials of degree less than or equal to $1$. In addition, the boundary conditions ($\A\p\bnu =0$ on $\partial\Omega$) are taken into account via penalization of the vortices that belong on the boundaries.Again thanks to geometric property of our computational domain, imposing this kind of boundary condition becomes trivial.   

	As stated earlier, the scan procedure consists in moving the probe according to the z-direction all the way through the tube. This leads to a set of scan positions where at each position we need to solve an ill-conditioned huge system. In addition to that, for any position one has to rebuild the triangulation and reassemble the matrices. To avoid this handicap we are led to consider a unique triangulation that includes all probe positions and then consider a factorization using a sophisticated software. 
	
\begin{figure}[!htbp]
\includegraphics[height=5cm,width=4cm]{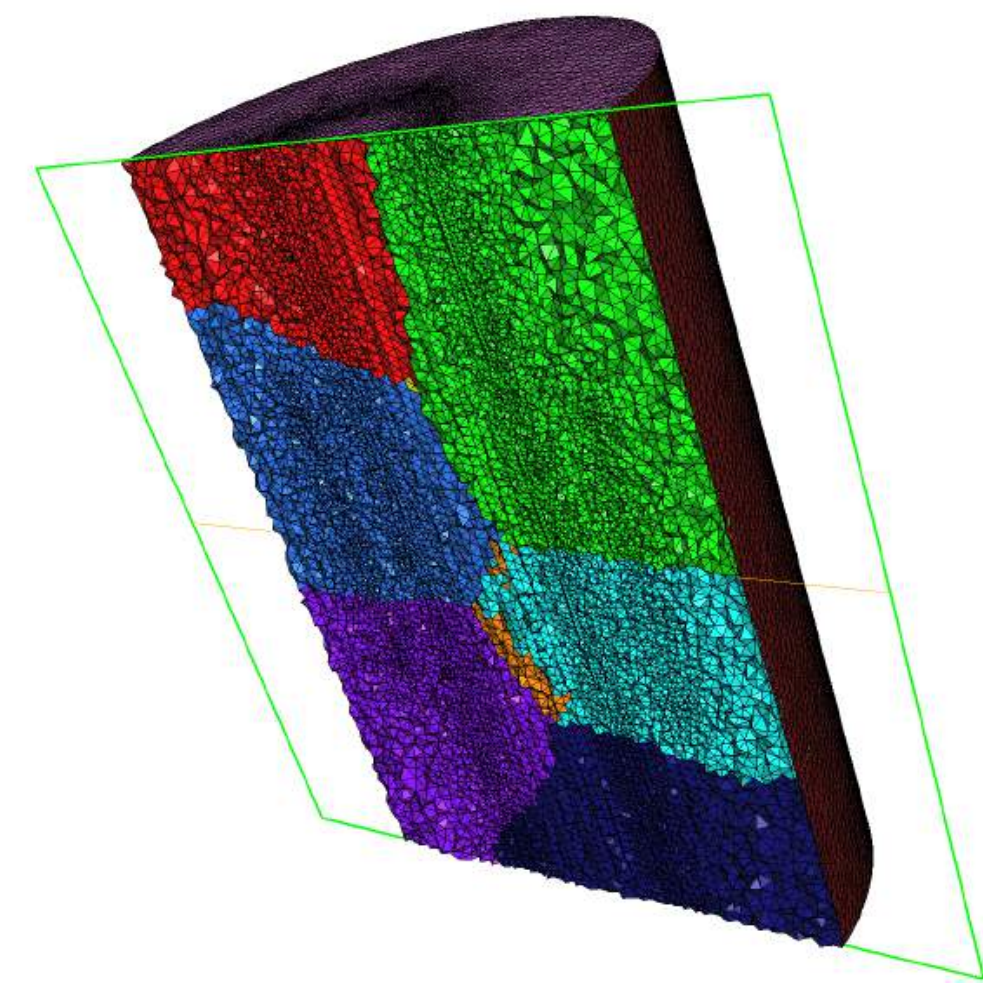}
\includegraphics[height=5cm,width=4cm]{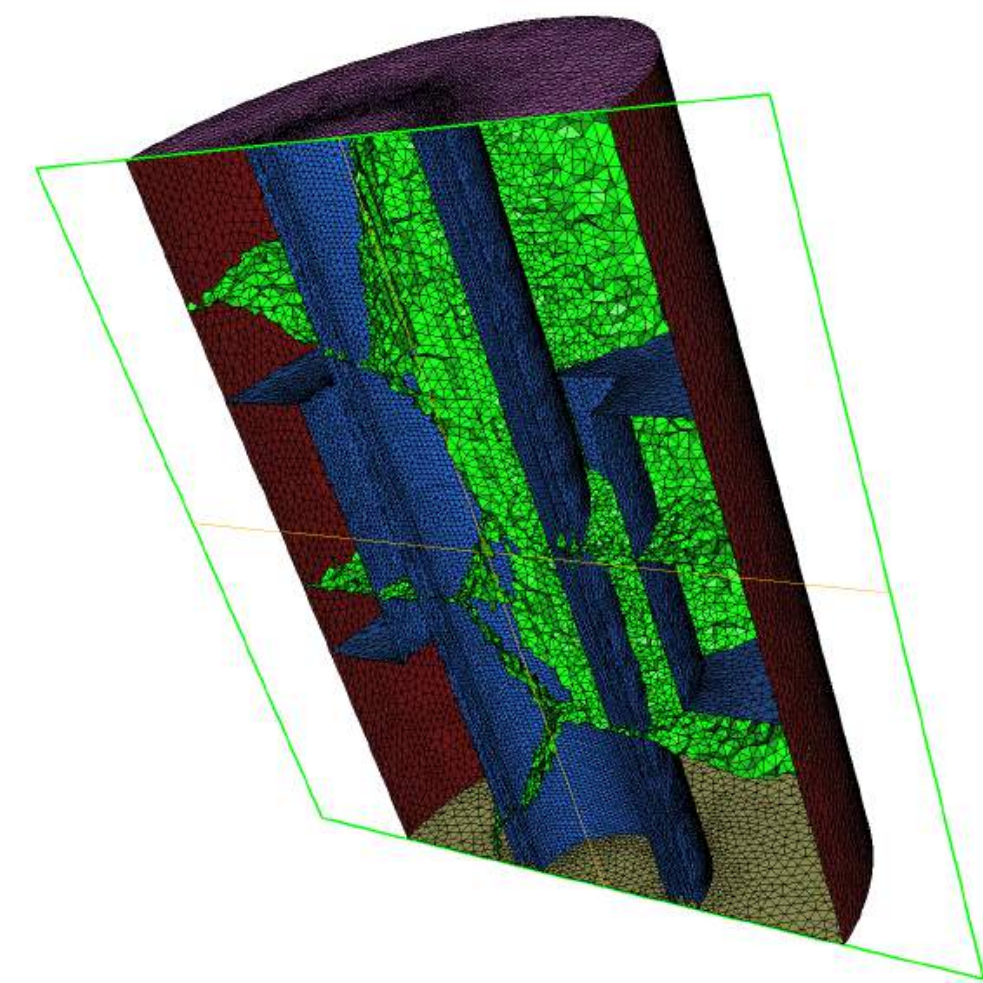}
\includegraphics[height=5cm,width=4cm]{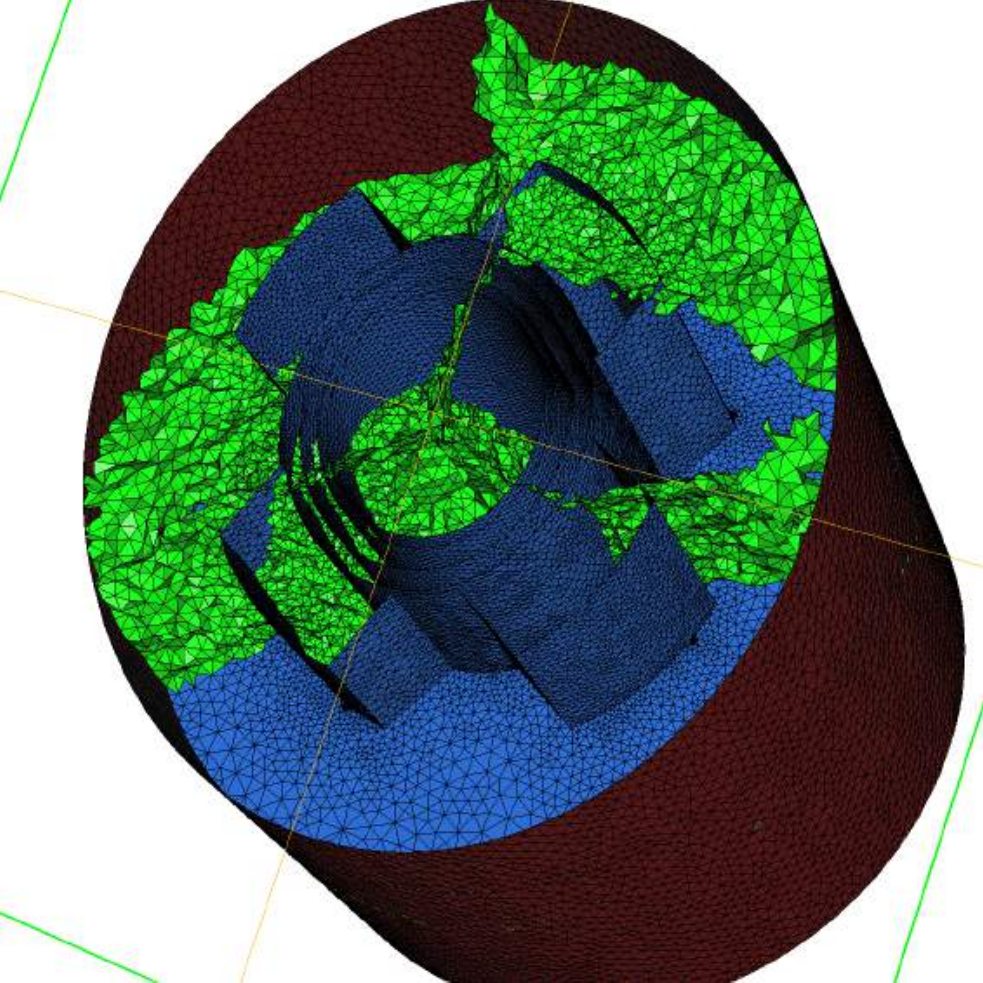}
\includegraphics[height=5cm,width=4cm]{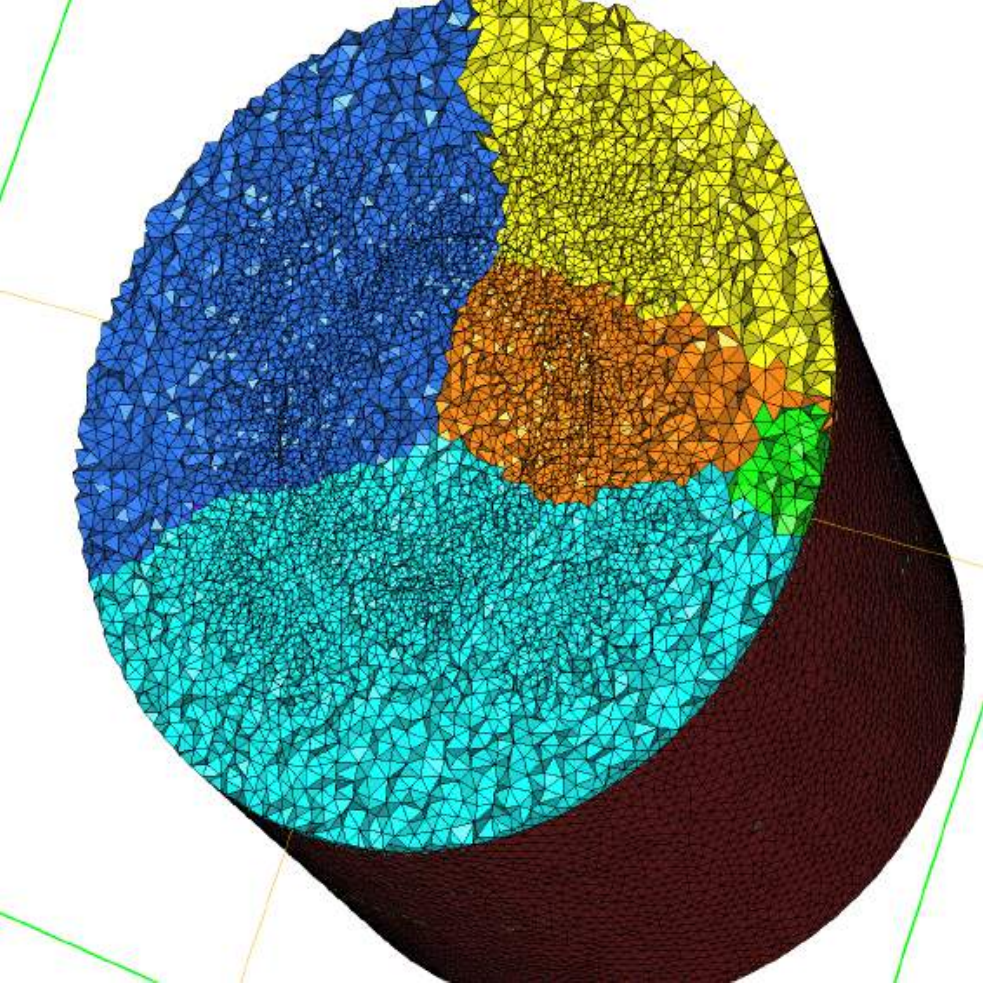}
\caption{Partition of the Triangulation of the computational domain. Volume tetrahedrons surfaces are presented in top left and bottom right, the other plot correspond to the interfaces as required by the geometry of the SG and also interfaces generated by the partition.}
\label{partitionTriangle}
\end{figure}

	In order to proceed with the parallelization, let us consider $P$ partitions of the computational domain $\Om$, such that $\Om=\cup_{p=1}^{P}\Om^p$. We notice here that the partition is made in terms of the triangulation see Figure~\ref{partitionTriangle}, where after the construction of the mesh, we partition its triangulations using a graph partitioner (e.g. Scotch~\cite{Pellegrini01scotchand} or Metis~\cite{Karypis95metis}). 
	\begin{remark}
	Both graph partitioners Scotch and Metis are interfaced with the software Freefem++~\cite{MR3043640}, and their utilization is quite easy. 
	\end{remark}
We give a sample Freefem++ script in Table~\ref{partition}, where the graph partition scotch is called at line 6 to build a piecewise function "balance" which takes different values on each subdomain. Freefem++ is therefore able to construct local triangulation according to the values of the function "balance" (see lines 11-12 of Table~\ref{partition}).  We notice that the construction of the function "balance" is performed only on the master processor with rank zero. The function "balance" is thus broadcasted to all processors through MPI (see line 9 of Table~\ref{partition}). We use this technique to ensure that all processors receive the same function balance and construct a uniform triangulation. 
	
	The numerical simulation of the ECT procedure reads as follows
\begin{enumerate}[Step 1.] \setlength\itemsep{-.6mm}
\item Read/build a triangulation of the computational domain $\Om$ 
\item Partition the triangulation with respect to the arbitrary choice (see command line in Tabular~\ref{partition})
\item Assemble Morse type sparse-sub-matrices in parallel
\item All\_reduce sparse-sub-matrices (This operation is achieved with MPI\_SUM operation) 
\item LU factorization of the global matrix.
\item Solve the set of problems $P_i$ by uniquely changing the right-hand side source term (accordingly to the position of the probe position)
\item Print out the impedance value results in an append writing file (to be sorted). 
\end{enumerate}

Let $\Phi_\natural\in \mathcal{V}_{\natural}$ and $\varphi_\natural\in \mathcal{C}_{\natural}$ be test functions for our coupled eddy-current problem, we denote by $\A_{\natural}$ and $\Vh$ our unknowns. Local matrices need the following variational formulation in order to be assembled
\begin{eqnarray} 
\mathcal{L}^{p}_{11}\big(\A_{\natural},\Phi_{\natural}\big)\!\!\!\!\!\!&=&\!\!\!\!\!\!\sum_{K\in\mathcal{T}_h\subset\Om^{p}}\int_{K}\dfrac{1}{\mu}\curl\A_{\natural}\p\overline{\curl\Phi_{\natural}} \dv,\notag\\
&\qquad&+\dfrac{1}{\tilde\mu}\int_{K}\dvg\A_{\natural}\overline{\dvg\Phi_{\natural}} \dv \\
\mathcal{L}^{p}_{12}\big(\Vh,\Phi_{\natural}\big)\!\!\!\!\!\!&=&\!\!\!\!\!\!-\!\!\!\!\sum_{K\in\mathcal{T}_h\subset\Om_c^{p}}\int_{K}\sigma \Vh\p\overline{\Phi_{\natural}} \dv.\\
\mathcal{L}^{p}_{21}\big(\A_{\natural},\varphi_{\natural}\big)\!\!\!\!&=&\!\!\!\!-\!\!\!\!\sum_{K\in\mathcal{T}_h\subset\Om_c^{p}}\int_{K}\sigma\A_{\natural}\p \overline{\grad\varphi_{\natural}} \dv,\\
\mathcal{L}^{p}_{22}\big(\Vh,\varphi_{\natural}\big)\!\!\!\!\!\!&=&\!\!\!\!\!\!-\dfrac{1}{i\omega}\!\!\!\!\!\!\sum_{K\in\mathcal{T}_h\subset\Om_c^{p}}\int_{K}\sigma\Vh \p\overline{\grad\varphi_{\natural}} \dv.
\end{eqnarray}
We give in Table~\ref{varf} a sample-freefem script that shows how to implement such bilinear forms.

\begin{equation*}
\begin{array}{lr}
 {\bf M}_{11} = \sum_{p=1}^{P}M^{p}_{11}, &\text{ with } { M}^{p}_{11}= \mathcal{L}^{p}_{11}(\Phi_{\natural},\Phi_{\natural}),\\
 {\bf M}_{12} = \sum_{p=1}^{P}M^{p}_{12}, &\text{ with } { M}^{p}_{12}=\mathcal{L}^{p}_{12}(\varphi_{\natural},\Phi_{\natural}),\\
 {\bf M}_{21} = \sum_{p=1}^{P}M^{p}_{21}, &\text{ with } { M}^{p}_{21}=\mathcal{L}^{p}_{21}(\Phi_{\natural},\varphi_{\natural}),\\
 {\bf M}_{22} = \sum_{p=1}^{P}M^{p}_{22}, &\text{ with } { M}^{p}_{22}=\mathcal{L}^{p}_{22}(\varphi_{\natural},\varphi_{\natural}).
\end{array}
\end{equation*}
\begin{remark}
The sparse matrice summations above are performed through the MPI Reduce operation. Thanks to the Morse sparse format of the matrices, the summations is not a term-by-term addition, but it increases the dimension at each sum, because of the assembly in disjoint subdomains; also all assembly uses the same mesh numbering. 
\end{remark}

The full discretized system thus reads
\begin{equation*}
\left(\begin{array}{cc}
{\bf M}^{}_{11} & {\bf M}^{}_{12}\\
{\bf M}{}_{21} & {\bf M}^{}_{22}
\end{array}\right)
\left(\begin{array}{l}\A_{\natural}\\\Vh\end{array}\right) =\left(\begin{array}{l} {\J}_{,\natural}\\ {\bf 0}\end{array}\right) 
\end{equation*}

The right hand side vector ${\J}_{,\natural}$ stands for the source term, which is basically given by the multiplication of the mass matrix by the interpolation of the analytic source term $\J$.

\begin{table}[!htbp]
\begin{lstlisting}
 int[int] nupart(Th.nt);
 fespace p0h(Th,P03d); 
 p0h balance;
 int npart=mpisize;
 if(mpirank==0){ 
 	   scotch(nupart, Th, npart);
	 for(int i=0;i<Th.nt;i++)
 		balance[][i] = nupart[i];}
 broadcast(processor(0,com),balance[]);
 Th3[int] Thpart(npart);
 for(int i=0;i<npart;i++) {
   Thpart[i]=trunc(Th,balance==i);
   Thpart[i]=change(Thpart[i],fregion=i);}
 Th=Thpart[0];
 for(int i=1;i<npart;i++) Th=Th+Thpart[i];
\end{lstlisting}
\caption{Mesh partitioning -- Freefem++ script sample.}\label{partition}
\end{table}

\begin{table}[!htbp]
\begin{lstlisting}
varf L11([Ax,Ay,Az],[Bx,By,Bz]) =
   int3d(Th,mpirank)(1/mu*curl(Bx,By,Bz)*curl(Ax,Ay,Az) )
    - int3d(Th,mpirank)( iomega*mu*sigma* [Bx,By,Bz]*[Ax,Ay,Az] )
    + int3d(Th,mpirank)((div(Bx,By,Bz))*(div(Ax,Ay,Az)))
    + on(labelup,labeldown,Az=0.)
    + on(labelmid,Ax=0., Ay=0.);
varf L22(V,qhc) =
    int3d(ThC,mpirank)( (-1/iomega) * (
     mu*sigmaEpsilon* [dx(qhc),dy(qhc),dz(qhc)]*[dx(V),dy(V),dz(V)] )
	+ delta*mu*sigma*V*qhc ) );
varf L12([V],[Bx,By,Bz]) =
    - int3d(ThC,mpirank)( mu*sigma* ( dx(V)*Bx+dy(V)*By+dz(V)*Bz) );
varf L21([Ax,Ay,Az],[qhc]) =
    - int3d(ThC,mpirank)( mu*sigma* [dx(qhc),dy(qhc),dz(qhc)]*[Ax,Ay,Az] );
\end{lstlisting}
\caption{Variational formulation definition for the coupled problem -- Freefem++ script sample.}\label{varf}
\end{table}

\begin{table}[!htbp]
\begin{lstlisting}
matrix <complex> M11,Mp11= L11(VPh,VPh); 
matrix <complex> M22,Mp22= L22(PhC,PhC);
matrix <complex> M12,Mp12= L12(PhCS,VPh);
matrix <complex> M21,Mp21= L21(VPh,PhCS);
\end{lstlisting}
\caption{Matrix Assembly, with Morse sparse forma -- Freefem++ script sample.}
\end{table}

\begin{table}[!htbp]
\begin{lstlisting}
mpiAllReduce(Mp11,M11,mpiCommWorld,mpiSUM);
mpiAllReduce(Mp12,M12,mpiCommWorld,mpiSUM);
mpiAllReduce(Mp21,M21,mpiCommWorld,mpiSUM);
mpiAllReduce(Mp22,M22,mpiCommWorld,mpiSUM);
matrix<complex> M = [[M11,M12],
                     [M21,M22]];
set(M,solver=sparsesolver,eps=1.e-16); 
\end{lstlisting}
\caption{Reduce Morse sparse matrix -- Freefem++ script sample.}
\end{table}

The boundary condition $\A\cdot\bnu=0$ on the exterior boundary $\partial\Om$ are taken into account using the penalization term, where the linear system is forced to take into account the given values at the boundary mesh points. 
\section{\normalsize Numerical Experiments}
	 In our application we consider that the excitation source term $\J$ is uniformly distributed on a support included in $\Om_I$ (principally it models the solenoid source coil of the probing problem). So, let $\J\in (L^2(\Om))^3$  with $\dvg\J=0$ in $\Om_I$ as $\J:=:{[-y_{|_{\Om_I}},x_{|_{\Om_I}},0]}\slash{\sqrt{x^2+y^2}}$. 
 
	As the electric scalar potential $\Vv$ is determined up to an additive constant, we may numerically impose a supplement condition such that
$\int_{\Om_{c_i}} \!\!\!\!\!\!\V \dv = 0$ ($\Om_{c_i}$ is any connex component subset of $\Om_c.$). This also could be incorporated under the global problem by penalization $\delta \sigma\V, \text{ in } \Om_c, \text{ with a small } \delta<<1.$ Therefore we augment the variational formulation $\mathcal{L}^{p}_{22}$ by the former penalization integration. 

\begin{table}[htbp]
  \centering
  \begin{tabular}{lcccc}
    \hline
    P in MPI& $\text{Assemble}\left(M^0_i\right)$ & $\text{Reduce}\left(M^0\right)$ & $LU|^0$\\    \hline
    1   & 3144.68    & -- & 585.28 \\
    2   &  1202.98   & 585.74 & 440.36\\
    4   & 758.73      & 575.86 & 582.76 \\
    8   &  312.08     & 512.76 & 501.25 \\
    16  & 163.87     & 591.56 &  596.33  \\
    32  & 88.74       & 542.30 & 577.28\\
    64  & 19.02 & 516.50 &  562.61\\
    \hline
  \end{tabular}
  \caption{Parallel Performance in wall-clock time in seconds for the problem without TSP as default. This means the numerical simulation takes the region of the TSP as a "vacuum".}  \label{TableM0}
\end{table}

\begin{table}[htbp]
  \centering
  \begin{tabular}{lcccc}
    \hline
    P in MPI& $\text{Assemble}\left(M^1_i\right)$ & $\text{Reduce}\left(M^1\right)$ & $LU|^1$  \\
    \hline
    1   &3237.47& - & 341.22  \\
        2   & 1284.52 & 392.18  & 332.61 \\
    4   & 764.47 & 587.49 & 332.44  \\
    8   & 311.28 & 687.63 &  497.75\\
    16   & 177.14 & 647.32 & 368.43 \\   
    32  &42.37 & 583.32 & 	328.61 &\\
    64  & 19.55 & 598.87 & 350.10 \\
    \hline
  \end{tabular}
    \caption{Parallel Performance in wall-clock time in seconds for the problem with TSP.}  \label{TableM1}
\end{table}

\begin{figure}[!htbp]
\centering
\includegraphics[height=4cm,width=4cm]{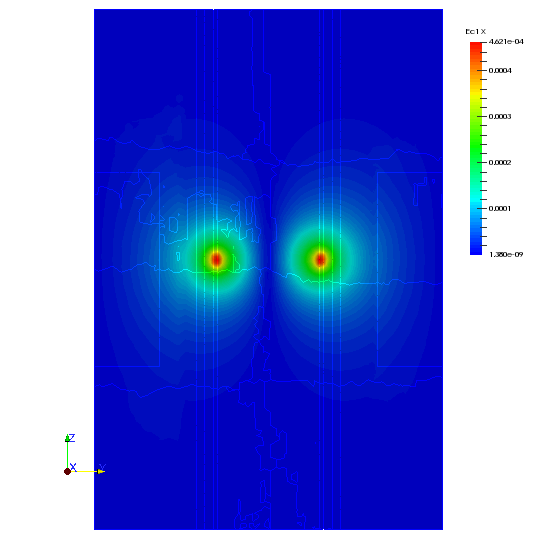}
\includegraphics[height=4cm,width=4cm]{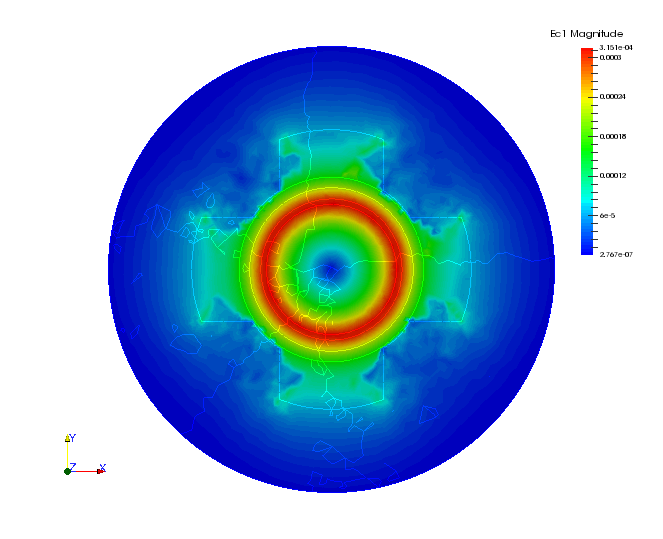}
\includegraphics[height=4cm,width=4cm]{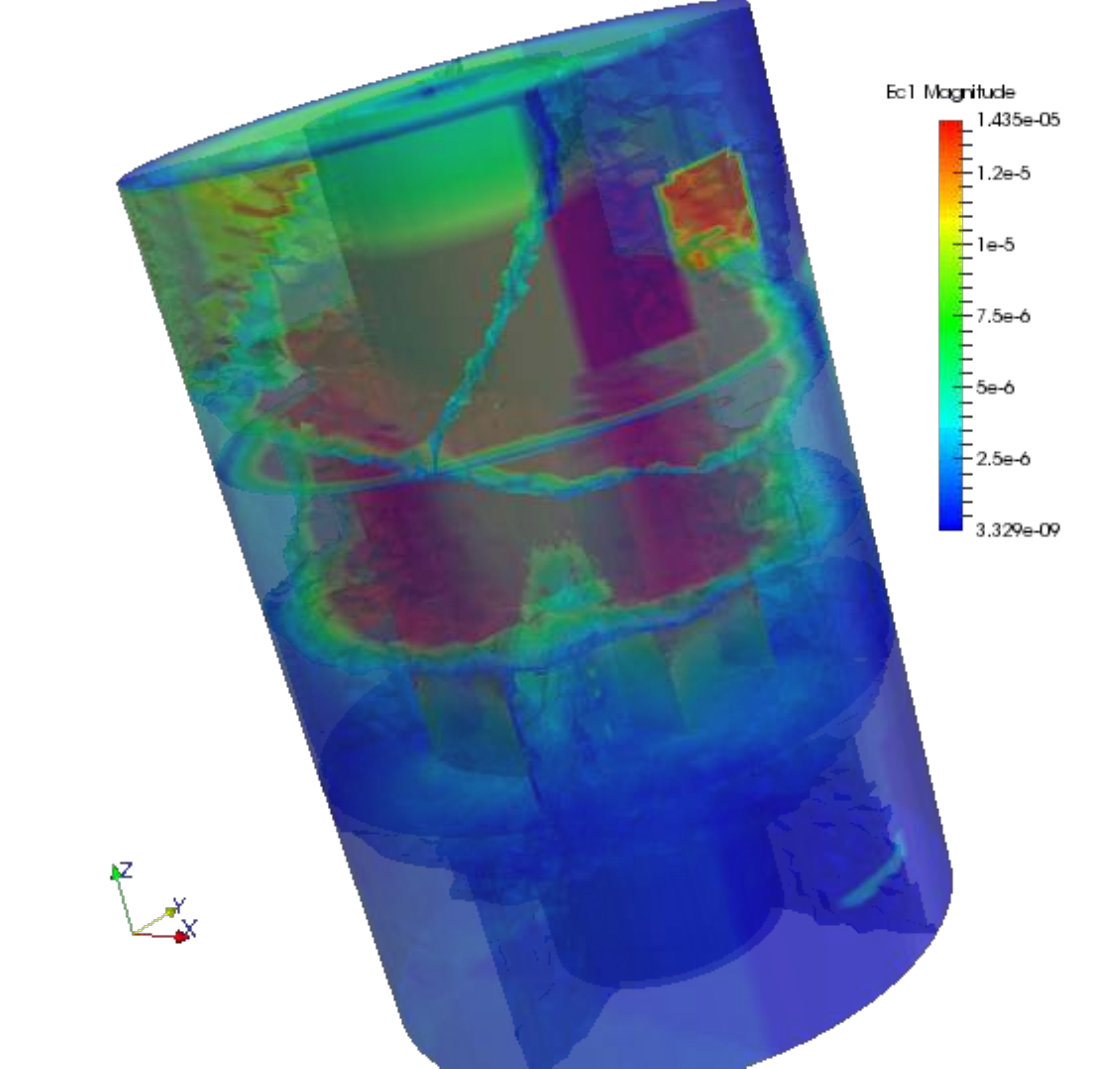}
\includegraphics[height=4cm,width=4cm]{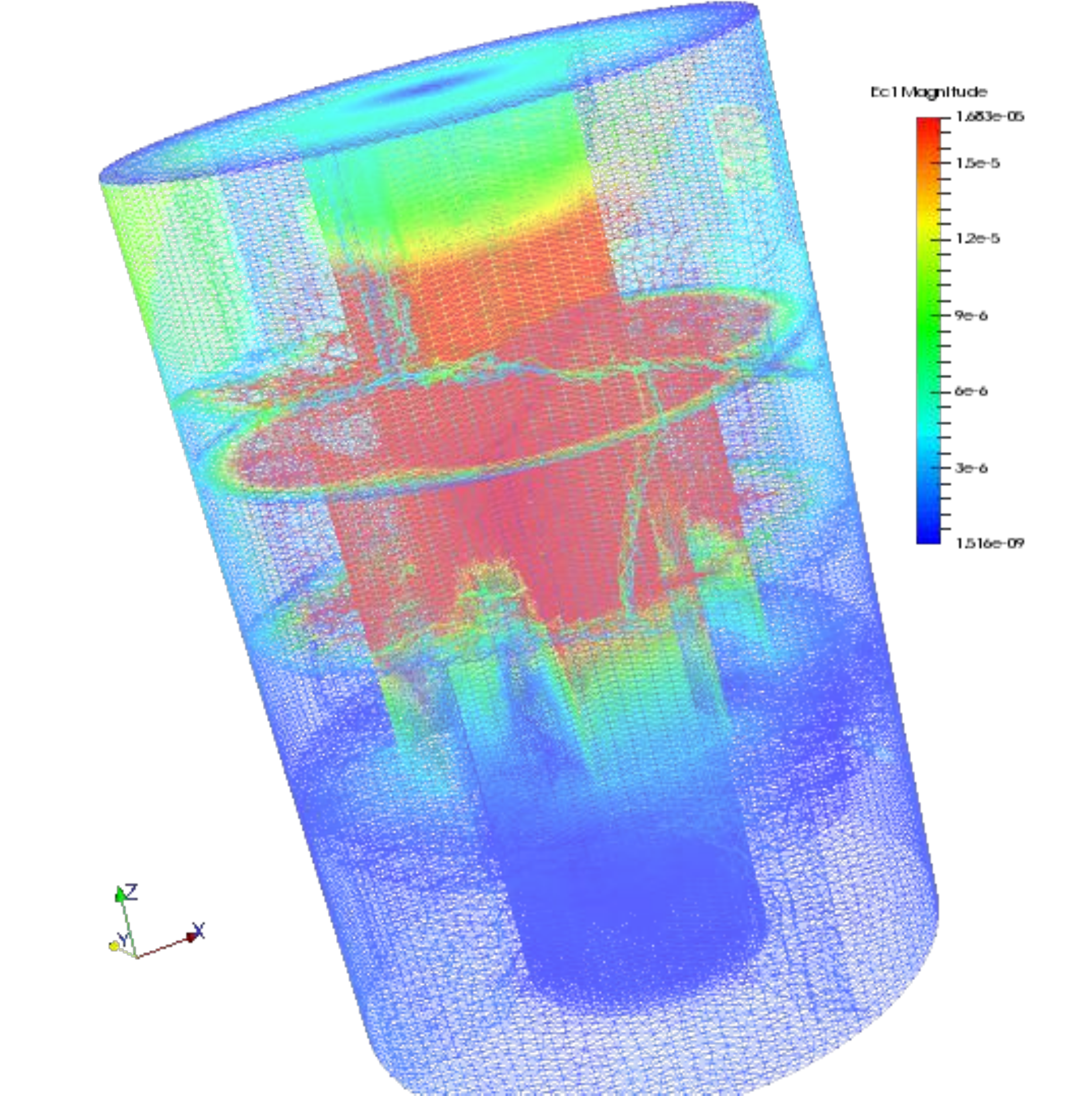}
\caption{The intensity of the computed Electric field $\E=i\omega\A+\V$: YZ-plan slice projection (top left) and XY-plan slice projections (top right). Volume rescale plot (bottom left) and Wireframe projection rescaled plot (bottom right).}\label{PlotE}
\end{figure}

\begin{figure}[!htbp]
\includegraphics[height=6cm,width=6cm]{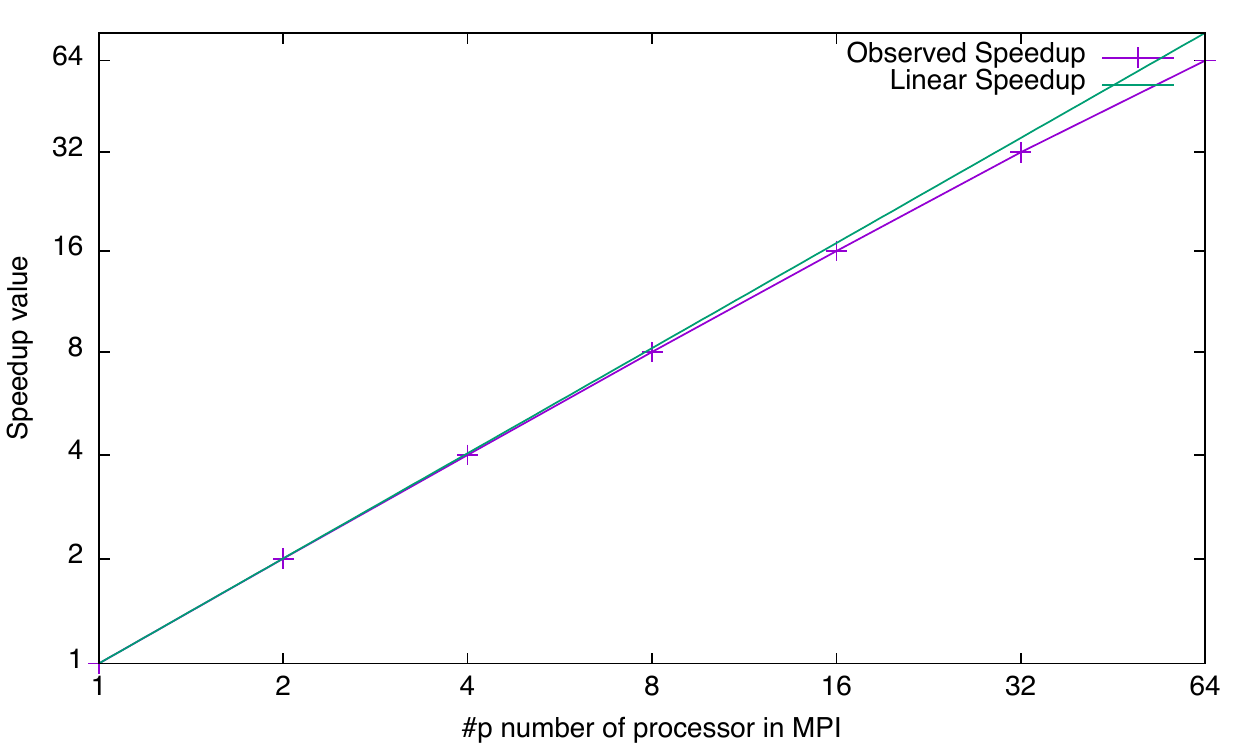}
\caption{Speedup of the numerical method.}\label{speedup}
\end{figure}

We give in Table~\ref{TableM0} and Table~\ref{TableM1} the performance in term of the execution (wall-clock) time in seconds of our numerical method. In Table~\ref{TableM0} (respectively Table~\ref{TableM1}) we report results related to the system without TSP as default (respectively with TSP region as default), where the computed matrix has size $n\times m=991.246\times 991.246$ (respectively $n\times m=1.069.595\times1.069.595$) with $46.664.492$ (respectively $54.350.718$) non-zero coefficients.
It is worth recalling that these calculations are necessary for the evaluation of the impedance signals Eq.~\eqref{impedkl}. As the results indicate, the most memory and time consuming task (the main matrix for the coupled system) is mitigated through high performance computing using parallel resources, where scalability of the operation is clearly exhibited. MPI communications are thus necessary to collect local Morse sparse matrices. These operations, thanks to the optimized communication technique, enable the maintenance of a reasonable average of performance despite the increase of the communicators number. Actually, this fact is balanced with the size of the message that has to be transferred. Indeed, as the communicators numbers increases, the size of the morse sparse matrices decreases. Once the global matrix is assembled and copied (via MPI) to all processors, we are then able to perform a fast factorization. In Figure~\ref{PlotE} we plot the intensity of the computed electric field when the probe is near to the TSP.
	We have used Super\_LU~\cite{li05} parallel direct solver for the factorization. This procedure has also a good performance while using several processors.   
As the scan procedure is parallel by nature (each probe position is totally independent of other positions). Therefore, parallel resources share the effort of solving independently coupled systems according to different probe positions. We report in Figure~\ref{speedup} the robustness of our method in term of speedup while increasing the number of parallel processors $P$.  We evaluate the speedup formula given by $S^{p}=\text{t}_\text{tot}^{p}\slash t^\text{serial}$ with $\text{t}_\text{tot}^{p}  = \text{t}_\text{Assemble}^{p} + \text{t}_\text{Reduce}^{p} + \text{t}_\text{Factorize}^{p} + \text{t}_\text{Solve}^{p}$. We have considered $128$ probe positions for the scan procedure where the time to solve the problem is in average $40$ minutes. The results exhibit the fact that our method is fully efficient and scalable in term of parallelism. 

\section{\normalsize Conclusions}\label{Conclusions}
We have presented in this paper a finite element technique that enables a rapid numerical simulation of the ECT problem in SG. Our approach has been validated through the use of high performance computing in parallel machines. We have shown that our approach is full efficient and leads to a robust and rapid ECT in real life industrial problem. 
\section*{\normalsize Acknowledgments}\label{Acknowledgements}
This work was made possible by the facilities of the New Jersey Institute of Technology Shared Hierarchical Academic Research Computing Network.
The author gratefully thanks Professor Frederic Hecht for the earliest discussion on Freefem++-mpi and also for the helpful introductory examples uploaded on the Freefem++ website.



\clearpage

\bibliographystyle{plainnat}
\bibliography{biblio}

\begin{thebibliography}{21}
\providecommand{\natexlab}[1]{#1}
\providecommand{\url}[1]{\texttt{#1}}
\expandafter\ifx\csname urlstyle\endcsname\relax
  \providecommand{\doi}[1]{doi: #1}\else
  \providecommand{\doi}{doi: \begingroup \urlstyle{rm}\Url}\fi

\bibitem[Alonso~Rodr{\'{\i}}guez and Valli(2010)]{MR2680968}
Ana Alonso~Rodr{\'{\i}}guez and Alberto Valli.
\newblock \emph{Eddy current approximation of {M}axwell equations}, volume~4 of
  \emph{MS\&A. Modeling, Simulation and Applications}.
\newblock Springer-Verlag Italia, Milan, 2010.
\newblock ISBN 978-88-470-1505-0.
\newblock \doi{10.1007/978-88-470-1506-7}.
\newblock URL \url{http://dx.doi.org/10.1007/978-88-470-1506-7}.
\newblock Theory, algorithms and applications.

\bibitem[Alonso~Rodr{\'i}guez et~al.(2013)Alonso~Rodr{\'i}guez, Bertolazzi,
  Ghiloni, and Valli]{MR3090156}
Ana Alonso~Rodr{\'i}guez, Enrico Bertolazzi, Riccardo Ghiloni, and Alberto
  Valli.
\newblock Construction of a finite element basis of the first de {R}ham
  cohomology group and numerical solution of 3{D} magnetostatic problems.
\newblock \emph{SIAM J. Numer. Anal.}, 51\penalty0 (4):\penalty0 2380--2402,
  2013.
\newblock ISSN 0036-1429.
\newblock \doi{10.1137/120890648}.
\newblock URL \url{http://dx.doi.org/10.1137/120890648}.

\bibitem[Ammari et~al.(2000)Ammari, Buffa, and
  N{\'e}d{\'e}lec]{Ammari:2000:JEC:354423.354457}
H.~Ammari, A.~Buffa, and J.-C N{\'e}d{\'e}lec.
\newblock A justification of eddy currents model for the maxwell equations.
\newblock \emph{SIAM J. Appl. Math.}, 60\penalty0 (5):\penalty0 1805--1823, May
  2000.
\newblock ISSN 0036-1399.
\newblock \doi{10.1137/S0036139998348979}.
\newblock URL \url{http://dx.doi.org/10.1137/S0036139998348979}.

\bibitem[Auld and Moulder(1999)]{auld1999review}
B.~A. Auld and J.~C. Moulder.
\newblock Review of advances in quantitative eddy current nondestructive
  evaluation.
\newblock \emph{Journal of Nondestructive evaluation}, 18\penalty0
  (1):\penalty0 3--36, 1999.

\bibitem[Bendjoudi et~al.(2010)Bendjoudi, Bossy, Cugnet, Chauvin, and
  Cassereau]{bendjoudi}
Aniss Bendjoudi, Emmanuel Bossy, Marie-Fran{\c c}oise Cugnet, Patrick Chauvin,
  and Didier Cassereau.
\newblock {D{\'e}veloppement d'un logiciel hybride pour le Contr{\^o}le Non
  Destructif}.
\newblock In Soci{\'e}t{\'e}~Fran{\c c}aise d'Acoustique SFA, editor,
  \emph{{10{\`e}me Congr{\`e}s Fran{\c c}ais d'Acoustique}}, pages~--, Lyon,
  France, 2010.
\newblock URL \url{http://hal.archives-ouvertes.fr/hal-00549199}.

\bibitem[B{\'{\i}}r{\'o} and Valli(2007)]{MR2298698}
Oszk{\'a}r B{\'{\i}}r{\'o} and Alberto Valli.
\newblock The {C}oulomb gauged vector potential formulation for the
  eddy-current problem in general geometry: well-posedness and numerical
  approximation.
\newblock \emph{Computer Methods in Applied Mechanics and Engineering},
  196\penalty0 (13-16):\penalty0 1890--1904, 2007.
\newblock ISSN 0045-7825.
\newblock \doi{10.1016/j.cma.2006.10.008}.
\newblock URL \url{http://dx.doi.org/10.1016/j.cma.2006.10.008}.

\bibitem[Biro and Valli(2007)]{biro2007coulomb}
Oszkar Biro and Alberto Valli.
\newblock The coulomb gauged vector potential formulation for the eddy-current
  problem in general geometry: well-posedness and numerical approximation.
\newblock \emph{Computer methods in applied mechanics and engineering},
  196\penalty0 (13):\penalty0 1890--1904, 2007.

\bibitem[Bossavit(2004)]{bossavit2004electromagnetisme}
Alain Bossavit.
\newblock \emph{Electromagnetisme, en vue de la modelisation}, volume~14.
\newblock Springer Science \& Business Media, 2004.

\bibitem[Durufl{\'e} et~al.(2006)Durufl{\'e}, Haddar, and
  Joly]{durufle2006higher}
Marc Durufl{\'e}, Houssem Haddar, and Patrick Joly.
\newblock Higher order generalized impedance boundary conditions in
  electromagnetic scattering problems.
\newblock \emph{Comptes Rendus Physique}, 7\penalty0 (5):\penalty0 533--542,
  2006.

\bibitem[Girard(2014)]{girarclogging}
Sylvain Girard.
\newblock Clogging of recirculating nuclear steam generators.
\newblock In \emph{Physical and Statistical Models for Steam Generator Clogging
  Diagnosis}, pages 3--13. Springer, 2014.

\bibitem[Haddar and Jiang(2015)]{haddar2015axisymmetric}
Houssem Haddar and Zixian Jiang.
\newblock Axisymmetric eddy current inspection of highly conducting thin layers
  via asymptotic models.
\newblock \emph{Inverse Problems}, 31\penalty0 (11):\penalty0 115005, 2015.

\bibitem[Haddar and Riahi(2013)]{haddar:hal-01044648}
Houssem Haddar and Mohamed~Kamel Riahi.
\newblock {3D direct and inverse solvers for eddy current testing of deposits
  in steam generator}.
\newblock Technical report, July 2013.
\newblock URL \url{https://hal.inria.fr/hal-01044648}.

\bibitem[Hecht(2012)]{MR3043640}
F.~Hecht.
\newblock New development in freefem++.
\newblock \emph{Journal of Numerical Mathematics}, 20\penalty0 (3-4):\penalty0
  251--265, 2012.
\newblock ISSN 1570-2820.

\bibitem[Huang and Takagi(2000)]{HuangTakagi}
Haoyu Huang and Toshiyuki Takagi.
\newblock Crack shape reconstruction from noisy signals in ect of steam
  generator tube.
\newblock In \emph{Industrial Electronics Society, 2000. IECON 2000. 26th
  Annual Conference of the IEEE}, volume~4, pages 2507--2512. IEEE, 2000.

\bibitem[Huang et~al.(2000)Huang, Takagi, and Fukutomi]{HuangTakagiFukutomi}
Haoyu Huang, T.~Takagi, and H.~Fukutomi.
\newblock Fast signal predictions of noised signals in eddy current testing.
\newblock \emph{Magnetics, IEEE Transactions on}, 36\penalty0 (4):\penalty0
  1719--1723, Jul 2000.
\newblock ISSN 0018-9464.
\newblock \doi{10.1109/20.877774}.

\bibitem[Karypis and Kumar(1995)]{Karypis95metis}
George Karypis and Vipin Kumar.
\newblock Metis - unstructured graph partitioning and sparse matrix ordering
  system, version 2.0, 1995.

\bibitem[Li(2005)]{li05}
Xiaoye~S. Li.
\newblock An overview of {SuperLU}: Algorithms, implementation, and user
  interface.
\newblock \emph{ACM Transactions on Mathematical Software}, 31\penalty0
  (3):\penalty0 302--325, September 2005.

\bibitem[Pellegrini(2001)]{Pellegrini01scotchand}
François Pellegrini.
\newblock Scotch and libscotch 3.4 user's guide, 2001.

\bibitem[Takagi et~al.(1997)Takagi, Tani, Fukutomi, and Hashimoto]{RPQNE}
Toshiyuki Takagi, Junji Tani, Hiroyuki Fukutomi, and Mitsuo Hashimoto.
\newblock Finite element modeling of eddy current testing of steam generator
  tube with crack and deposit.
\newblock In DonaldO. Thompson and DaleE. Chimenti, editors, \emph{Review of
  Progress in Quantitative Nondestructive Evaluation}, volume~16 of
  \emph{Review of Progress in Quantitative Nondestructive Evaluation}, pages
  263--270. Springer US, 1997.
\newblock ISBN 978-1-4613-7725-2.
\newblock \doi{10.1007/978-1-4615-5947-4_34}.
\newblock URL \url{http://dx.doi.org/10.1007/978-1-4615-5947-4_34}.

\bibitem[{The Open Access NDT Database}(2014)]{ndtdatabase}
{The Open Access NDT Database}.
\newblock The web's largest database of nondestructive testing ({NDT})
  conference proceedings, articles, news, exhibition, forum and a professional
  network.
\newblock \emph{NDT Database and Journal of Nondestructive Testing - NDT,
  Ultrasonic Testing, X-Ray, Radiography, Eddy Current and All NDT Methods.},
  2014.
\newblock URL \url{http://www.ndt.net/index.php}.

\bibitem[Tian et~al.(2005)Tian, Sophian, Taylor, and Rudlin]{1386227}
Gui~Yun Tian, A.~Sophian, D.~Taylor, and J.~Rudlin.
\newblock Multiple sensors on pulsed eddy-current detection for 3-d subsurface
  crack assessment.
\newblock \emph{Sensors Journal, IEEE}, 5\penalty0 (1):\penalty0 90--96, Feb
  2005.
\newblock ISSN 1530-437X.
\newblock \doi{10.1109/JSEN.2004.839129}.

\end{thebibliography}

%

\end{document}